\documentclass[preprint,12pt]{elsarticle}

\usepackage[T1]{fontenc}
\usepackage{lmodern}
\usepackage{amsmath,amssymb,amsthm}
\usepackage[shortlabels]{enumitem}
\usepackage{microtype}
\usepackage[hidelinks]{hyperref}

\newtheorem{theorem}{Theorem}[section]
\newtheorem{lemma}[theorem]{Lemma}
\newtheorem{proposition}[theorem]{Proposition}
\newtheorem{corollary}[theorem]{Corollary}
\theoremstyle{definition}
\newtheorem{definition}[theorem]{Definition}
\theoremstyle{remark}
\newtheorem{remark}[theorem]{Remark}

\newcommand{\FLe}{\mathrm{FL}_{e}}
\newcommand{\Inv}{\mathfrak I}
\newcommand{\Sym}{\mathfrak S}
\newcommand{\Div}{\mathfrak D}
\newcommand{\Chain}{\mathsf c}
\newcommand{\Sem}{\mathsf{sl}}
\newcommand{\odd}{\mathsf o}
\newcommand{\eveni}{\mathsf e_{\mathrm i}}
\newcommand{\evenn}{\mathsf e_{\mathrm n}}
\newcommand{\oddeven}{\mathsf{oe}}
\newcommand{\Z}{\mathbb Z}
\newcommand{\Idp}{\operatorname{Id}^{+}}

\begin{document}

\begin{frontmatter}

\title{Amalgamation in classes of involutive commutative residuated lattices\tnoteref{funding}}

\author{S\'andor Jenei}
\address{Institute of Mathematics and Informatics, Eszterh\'azy K\'aroly Catholic University, Hungary}
\address{Institute of Mathematics and Informatics, University of P\'ecs, Hungary}
\tnotetext[funding]{The present scientific contribution was supported by the NKFI-K-146701 grant.}

\begin{abstract}
We study amalgamation in odd and even involutive commutative residuated chains through their categorical representation by bunches of linearly ordered abelian groups.
The representation separates three obstructions.
A discrete $\kappa_J$-layer imports the failure of amalgamation for discrete abelian ordered groups with normal embeddings.
Even when $\kappa_J$ is empty, distinguished layer subgroups obstruct amalgamation in the unrestricted idempotent-symmetric classes.
Moreover, purity of every induced layer embedding does not suffice when the two targets enlarge the positive-idempotent skeleton in different ways.
We then isolate a sufficient positive regime.
For every finite $n\geq1$, the idempotent-symmetric odd chains, the idempotent-symmetric even chains with idempotent falsum, and their union have the Amalgamation Property after further restriction to algebras with exactly $n$ positive idempotents and divisible canonical layer groups.
Fixed $n$ identifies the skeletons, divisibility makes the layer embeddings pure, ordinary abelian-group pushouts remain torsion-free, and a simultaneous order-extension theorem produces an amalgamating bunch.
Although failure of SAP already follows formally from failure of AP, we also give an independent one-layer witness to failure of SAP in the unrestricted idempotent-symmetric odd and even classes.
At the variety level, the essential counterexamples yield a general transfer criterion for failure of AP in semilinear varieties; in particular, AP fails for the varieties generated by the idempotent-symmetric odd and even chains, and for every semilinear variety containing the variety of odd semilinear involutive commutative residuated lattices.
Finally, although every fixed-$n$, layer-divisible chain class considered here has AP, among the varieties they generate only the odd one-layer variety has AP.
\end{abstract}


\end{frontmatter}

\section{Introduction}\label{sec:introduction}

Residuation is an order-theoretic form of adjunction and has long provided a common language for ordered algebra; see \cite{BlythJanowitz,GalatosEtAl}.
The systematic study of residuated lattices dates back to Ward and Dilworth \cite{WardDilworth,Dilworth}, whose work was motivated in part by the ideal theory of rings developed by Krull and others \cite{Krull}.
The framework now encompasses lattice-ordered groups, Boolean and Heyting algebras, MV- and BL-algebras, and many related ordered structures.
Its modern logical importance comes from the fact that varieties of pointed residuated lattices serve as algebraic semantics for substructural logics \cite{GalatosEtAl}.

The Amalgamation Property and its variants are closely connected with interpolation and Robinson-style properties of the corresponding consequence relations.
Amalgamation first appeared in group theory through amalgamated free products \cite{Schreier}; in ordered algebra it has since been studied for groups, lattice-ordered groups, and a broad range of residuated structures.
A systematic account of the connections among amalgamation, congruence extension, and interpolation in ordered algebras is given in \cite{MetcalfeMontagnaTsinakis}.

Much of the residuated-lattice literature concerns linear, semilinear, conic, or semiconic varieties, frequently under integrality, divisibility, or idempotence assumptions; for a recent overview of amalgamation in the semilinear setting, see \cite{FussnerSantschiSemilinear}.
The general semilinear and semilinear cancellative varieties fail AP \cite{GilFerezLeddaTsinakis}, while significant positive structure and amalgamation results are available in idempotent settings \cite{Chen,GilFerezJipsenMetcalfe}.
A particularly close positive comparison is provided by the variety of odd Sugihara monoids, which has the Strong Amalgamation Property \cite{GalatosRaftery}.
The present paper addresses a complementary regime: odd and even involutive commutative residuated chains, without assuming integrality, divisibility of the algebra, or idempotence of multiplication.
These structures are algebraic semantics for involutive substructural logics without weakening, and their layer-group representation makes it possible to identify the precise sources of amalgamation failure.

The main tool is the representation of an odd or even involutive $\FLe$-chain by a bunch of layer groups \cite{JeneiRep,JeneiCat,JeneiDense}.
The positive idempotents form a totally ordered skeleton.
Each skeleton point carries a linearly ordered abelian group.
The $\kappa_J$-points carry discrete groups, whereas the $\kappa_I$-points carry distinguished subgroups.
The groups form a direct system, and algebra embeddings induce embeddings of the skeletons and of the layer groups.
At a $\kappa_J$-point they preserve neighboring elements.
At a $\kappa_I$-point they preserve both the distinguished subgroup and its complement.

These features yield three genuinely different obstructions.

First, the discrete $\kappa_J$-layers import a known non-amalgamable V-formation of discrete abelian ordered groups with positive normal embeddings \cite[Lemma~12]{GlassPierce}.
This proves failure of the Amalgamation Property (AP) for the full odd and even chain classes.

Second, deleting all $\kappa_J$-points does not restore AP.
At a $\kappa_I$-layer, an alleged amalgam can be forced to place the same group element simultaneously inside and outside its distinguished subgroup.
Thus the unrestricted idempotent-symmetric odd and even classes still fail AP.

Third, purity of all layer embeddings does not control the skeleton.
When the two targets insert different intermediate idempotents, the transition maps into their distinguished subgroups can force another contradiction.
This shows that a same-skeleton assumption in the representation-side construction cannot be replaced merely by purity of the arrows.

The positive result works inside the idempotent-symmetric odd and even classes and uses two further restrictions, each imposed on an individual algebra rather than on a chosen amalgamation diagram.
Fix $n\geq1$.
We require exactly $n$ positive idempotents and require every canonical layer group to be divisible.
In a V-formation inside this class, the induced maps between the finite $n$-element skeletons are automatically isomorphisms.
After identifying the skeletons, divisibility makes every induced layer embedding pure.
The ordinary abelian-group pushouts are then torsion-free and contain both factors as ordered subgroups.
Their transition maps are supplied by the universal property.
A simultaneous order-extension theorem makes all pushout layers linearly ordered without changing the groups or transitions.
The resulting bunch is the required amalgam and still has divisible layers.

The negative theory has stronger variety-level consequences.
A one-layer reduction transfers the failure of strong amalgamation exhibited for abelian linearly ordered groups in \cite[Theorem~2]{CherriPowell} to the idempotent-symmetric odd and even chain classes.
Both the discrete-layer counterexample and the distinguished-subgroup counterexamples can be chosen essential on one leg.
Combining these witnesses with the semilinear transfer theorem \cite[Corollary~3.5]{FussnerSantschi}, based on the general results of \cite{FussnerMetcalfe}, yields a useful transfer principle: any semilinear variety containing one of the odd essential V-formations constructed here fails AP; the analogous conclusion holds for the even V-formations when all chains in the variety are odd or even.
As consequences, every semilinear variety containing the variety of odd semilinear involutive commutative residuated lattices fails AP, and so do the varieties generated by the idempotent-symmetric odd chains, the idempotent-symmetric even chains, or their union.

This also exposes a distinction between a non-equational chain class and the variety it generates.
For each fixed finite $n$, the layer-divisible chain classes treated in the positive theorem have AP.
Nevertheless, their generated varieties fail AP in every even case and in every odd case with at least two positive idempotents.
The unique exception is the variety generated by the odd one-layer class: it is term-equivalent to the variety of abelian lattice-ordered groups and has AP \cite[Corollaries~6.2 and~6.3]{FussnerSantschiSemilinear}.

The resulting picture corrects the tempting but invalid inference that eliminating the discrete layers alone should yield AP.
The distinguished subgroups and the geometry of the positive-idempotent skeleton are independent sources of obstruction, and the positive theorem must control both.

Section~\ref{sec:preliminaries} fixes terminology and records the representation data used below.
Section~\ref{sec:negative} proves the three chain-level obstructions and the failure of strong amalgamation.
Section~\ref{sec:group-tools} develops the group-theoretic pushout and order-extension tools.
Section~\ref{sec:positive} proves the fixed-idempotent positive theorem.
Section~\ref{sec:varieties} gives the semilinear transfer and determines the AP status of the varieties generated by the positive chain classes.

\section{Preliminaries}\label{sec:preliminaries}

\subsection{Algebraic setting}

A \emph{commutative residuated lattice} is an algebra
$$
 \mathbf L=(L,\wedge,\vee,\cdot,\mathbin{\to},t)
$$
such that $(L,\wedge,\vee)$ is a lattice, $(L,\cdot,t)$ is a commutative monoid, and multiplication and the residual satisfy
$$
 xy\leq z
 \quad\Longleftrightarrow\quad
 y\leq x\to z
 \qquad(x,y,z\in L),
$$
where $\leq$ is the lattice order.
An $\FLe$-algebra is an expansion
$$
 \mathbf X=(X,\wedge,\vee,\cdot,\mathbin{\to},t,f)
$$
of a commutative residuated lattice by a distinguished constant $f$.
The residual complement is $x'=x\to f$.
The algebra is \emph{involutive} if $x''=x$, and it is a \emph{chain} if its lattice order is total.
An involutive $\FLe$-algebra is \emph{odd} if $f=t$.
It is \emph{even} if, for every $x$,
$$
 x<t \quad\Longleftrightarrow\quad x\leq f.
$$
For chains this is equivalent to saying that $f$ is the lower cover of $t$.
In the even case, the falsum may be idempotent or non-idempotent.
An $\FLe$-algebra is \emph{semilinear} if it is a subdirect product of $\FLe$-chains.

Write
$$
 \Inv_{\odd}^{\Chain},\qquad
 \Inv_{\eveni}^{\Chain},\qquad
 \Inv_{\evenn}^{\Chain}
$$
for the classes of odd chains, even chains with idempotent falsum, and even chains with non-idempotent falsum, respectively.
Put
$$
 \Inv_{\oddeven}^{\Chain}
 =
 \Inv_{\odd}^{\Chain}
 \cup
 \Inv_{\eveni}^{\Chain}
 \cup
 \Inv_{\evenn}^{\Chain}.
$$

An odd or even involutive $\FLe$-algebra is \emph{idempotent-symmetric} if
$$
 x^2=x
 \quad\Longrightarrow\quad
 (x')^2=x'.
$$
Write $\Sym_{\odd}^{\Chain}$ and $\Sym_{\eveni}^{\Chain}$ for the idempotent-symmetric odd and even chain classes.
Idempotent-symmetry forces the falsum of an even algebra to be idempotent, because $t$ is idempotent and $t'=f$.

For an odd or even involutive chain $\mathbf X$, let
$$
 \Idp(\mathbf X)
 =
 \{u\in X:u\geq t,\ u^2=u\}.
$$
This, equipped with the ordering inherited from $\mathbf X$, is the \emph{positive-idempotent skeleton} of $\mathbf X$.

\subsection{Amalgamation notions}

\begin{definition}\label{def:ap}
Let $\mathcal K$ be a class of similar structures.
A \emph{V-formation} in $\mathcal K$ is a pair of embeddings
$$
 \mathbf A\xrightarrow{\iota_1}\mathbf B,
 \qquad
 \mathbf A\xrightarrow{\iota_2}\mathbf C
$$
with all three structures in $\mathcal K$.
An \emph{amalgam} is a structure $\mathbf D\in\mathcal K$ and embeddings $\mu_1:\mathbf B\to\mathbf D$ and $\mu_2:\mathbf C\to\mathbf D$ satisfying $\mu_1\iota_1=\mu_2\iota_2$.
The class has the \emph{Amalgamation Property} if every V-formation in it has an amalgam in it.
\end{definition}

\begin{definition}\label{def:sap}
An amalgam is \emph{strong} if
$$
 \mu_1(B)\cap\mu_2(C)
 =
 \mu_1\iota_1(A)
 =
 \mu_2\iota_2(A).
$$
A class has the \emph{Strong Amalgamation Property} (SAP) if every V-formation in it has a strong amalgam in it.
\end{definition}

\subsection{Ordered groups and direct systems}

We write abelian groups additively whenever group calculations are involved.
A \emph{partially ordered abelian group} is an abelian group equipped with a translation-invariant partial order.
Its positive cone is
$$
 G^+=\{g\in G:0\leq g\}.
$$
A homomorphism $\varphi:G\to H$ is \emph{positive} if $\varphi(G^+)\subseteq H^+$, or equivalently if it is order-preserving.
An injective homomorphism is an \emph{ordered-group embedding} when it both preserves and reflects the order.
A partially ordered abelian group is \emph{linearly ordered} when its order is total; every such group is torsion-free.
A linearly ordered abelian group is \emph{discrete} if it has a least strictly positive element; equivalently, every element has an immediate predecessor and an immediate successor.

Let $\kappa$ be a directed poset.
A \emph{direct system} of abelian groups over $\kappa$ is a family
$$
 \bigl(G_u,\sigma_{uv}\bigr)_{u\leq v\text{ in }\kappa}
$$
such that
$$
 \sigma_{uu}=\operatorname{id}_{G_u},
 \qquad
 \sigma_{uw}=\sigma_{vw}\sigma_{uv}
 \quad(u\leq v\leq w).
$$
For two systems over the same index set, a morphism
$$
 \mathbf f:(G_u,\sigma^G_{uv})\longrightarrow(K_u,\sigma^K_{uv})
$$
is a family of homomorphisms $f_u:G_u\to K_u$ satisfying
$$
 f_v\sigma^G_{uv}=\sigma^K_{uv}f_u
 \qquad(u\leq v).
$$
When the groups are ordered, the transitions and component maps are required to be positive.
All skeletons below are chains and hence are directed.

\subsection{The representation package}

We recall only the features of the layer-group representation used in the proofs.
The complete construction and its categorical form are given in \cite{JeneiRep,JeneiCat,JeneiDense}.

A bunch has a chain $\kappa$ with least element $t$, partitioned as
$$
 \kappa
 =
 \kappa_{\odd}\mathbin{\dot\cup}\kappa_J\mathbin{\dot\cup}\kappa_I,
 \qquad
 \kappa_{\odd}\subseteq\{t\}.
$$
Each $u\in\kappa$ carries a linearly ordered abelian group $G_u$.
For $u\leq v$, there is a positive homomorphism
$$
\sigma_{uv}:G_u\to G_v,
$$
and these maps form a direct system.
For a discrete linearly ordered group, write $a_{\downarrow}$ and $a_{\uparrow}$ for the lower and upper neighbors of $a$.
Every $\kappa_J$-layer group is discrete, and whenever $u\in\kappa_J$ and $u<v$, one has
$$
 \sigma_{uv}(0)=\sigma_{uv}(0_{\downarrow}).
$$
Every $\kappa_I$-layer carries a distinguished subgroup $H_u\leq G_u$, and
$$
 \sigma_{vu}(G_v)\subseteq H_u
 \qquad(v<u,\ u\in\kappa_I).
$$

As part of the bunch data, the neutral element of $G_u$ is the skeleton point $u$, and the universes of all $G_u$'s and all dotted copies $\dot H_u$ introduced below are pairwise disjoint.
In additive calculations below we suppress the tags, write each neutral element as $0$, and understand every displayed group and dotted copy as a tagged disjoint isomorphic copy.

An ordered-group embedding between discrete groups is \emph{normal} if it preserves these neighbors.
Bunch embeddings induce normal embeddings on $\kappa_J$-layers.

At an $\kappa_I$-layer, a bunch embedding $F$ satisfies the exact inverse-image condition
\begin{equation}\label{eq:H-inverse-image}
 F_u^{-1}(H_{\alpha(u)})=H_u.
\end{equation}
Thus it preserves both the distinguished subgroup and its complement.

In the chain reconstructed from a bunch, the $\kappa_I$-layer corresponding to $u$ has universe
$$
G_u\mathbin{\dot\cup}\dot H_u,
\qquad
\dot H_u=\{\dot h:h\in H_u\},
$$
where $\dot H_u$ is a disjoint copy of $H_u$.
The elements of $G_u$ are called \emph{undotted}, and those of $\dot H_u$ are called \emph{dotted}.
Let $u\in\kappa_I$, and let $F_u:G_u\to G_{\alpha(u)}$ be the corresponding layer component of a bunch homomorphism.
If $\alpha(u)\in\kappa_I$, then the induced map on the reconstructed $\kappa_I$-layer is given by
$$
g\longmapsto F_u(g)\quad(g\in G_u),
\qquad
\dot h\longmapsto\dot{F_u(h)}\quad(h\in H_u).
$$
If $\alpha(u)\in\kappa_{\odd}$, then it is given by
$$
g\longmapsto F_u(g)\quad(g\in G_u),
\qquad
\dot h\longmapsto F_u(h)\quad(h\in H_u).
$$
For an embedding, the case $\alpha(u)\in\kappa_{\odd}$ cannot occur, since it would identify $h$ and $\dot h$.
Hence embeddings preserve the distinction between dotted and undotted elements.

We shall also use the following features of arbitrary bunch homomorphisms.
Their skeleton maps are isotone and preserve the least skeleton point.
A $\kappa_J$-point may be mapped only to a $\kappa_J$-point or to the unique $\kappa_{\odd}$-point, and a $\kappa_I$-point may be mapped only to a $\kappa_I$-point or to the unique $\kappa_{\odd}$-point.
If a $\kappa_J$-point is mapped to a $\kappa_J$-point, the corresponding layer homomorphism preserves the neighboring-element operations.
Finally, the partial-injectivity condition implies that if $u\in\kappa_I$, $\alpha(u)\in\kappa_I$, and $H_u=G_u$, then $F_u$ is strictly order-preserving and hence injective.
These facts are immediate special cases of the homomorphism conditions in the categorical representation cited below.

\begin{theorem}[Layer-group representation]\label{thm:representation}
The category of odd or even involutive $\FLe$-chains and their homomorphisms is categorically equivalent to the category of bunches of layer groups and their homomorphisms.
Under this equivalence:
\begin{enumerate}[(i)]
\item the skeleton of the bunch of $\mathbf X$ is
$$
 \kappa(\mathbf X)
 =
 \{x\to x:x\in X\}
 =
 \Idp(\mathbf X);
$$
\item $\mathbf X$ is odd precisely when $t\in\kappa_{\odd}$, even with idempotent falsum precisely when $t\in\kappa_I$, and even with non-idempotent falsum precisely when $t\in\kappa_J$;
\item $\mathbf X$ is idempotent-symmetric precisely when $\kappa_J=\varnothing$;
\item algebra embeddings correspond to bunch embeddings.
\end{enumerate}
\end{theorem}

\begin{proof}
The object-level representation is a reformulation of \cite[Theorem~8.1]{JeneiRep} in the notation and conventions used here.
The categorical correspondence is proved in \cite[Theorem~3.6]{JeneiCat}.
The embedding characterization used here is recorded in \cite[Lemma~10]{JeneiDense}.
\end{proof}

For later use, we isolate the $\kappa_J=\varnothing$ data.

\begin{definition}\label{def:symmetric-bunch}
A \emph{symmetric bunch} consists of:
\begin{enumerate}[(B1),leftmargin=3em]
\item a chain $\kappa$ with least element $t$, partitioned as
$$
 \kappa=\kappa_{\odd}\mathbin{\dot\cup}\kappa_I,
 \qquad
 \kappa_{\odd}\subseteq\{t\};
$$
\item a linearly ordered abelian group $G_u$ for every $u\in\kappa$;
\item positive homomorphisms
$$
 \sigma_{uv}:G_u\to G_v
 \qquad(u\leq v),
$$
satisfying the direct-system identities;
\item a subgroup $H_v\leq G_v$ for every $v\in\kappa_I$, such that
$$
 \sigma_{uv}(G_u)\subseteq H_v
 \qquad(u<v).
$$
\end{enumerate}
A symmetric-bunch embedding consists of an order embedding $\alpha$ of the skeletons satisfying
$$
 \alpha(t)=t,\qquad
 \alpha(\kappa_{\odd})\subseteq\kappa_{\odd},
 \qquad
 \alpha(\kappa_I)\subseteq\kappa_I,
$$
together with compatible ordered-group embeddings on the layers satisfying \eqref{eq:H-inverse-image}.
\end{definition}

\section{Negative chain-level results}\label{sec:negative}

\subsection{\texorpdfstring{The discrete $\kappa_J$-layer obstruction}{The discrete kappa-J-layer obstruction}}

The first obstruction comes from discrete abelian ordered groups.
We use the following known result.

\begin{lemma}\label{lem:discrete-group-obstruction}
There exists a V-formation of discrete abelian linearly ordered groups and positive normal embeddings which has no amalgam by positive normal embeddings in a discrete abelian linearly ordered group.
\end{lemma}

\begin{proof}
Glass and Pierce work with discrete abelian linearly ordered groups in a language containing a distinguished constant $1$, axiomatized as the least strictly positive element, and prove that this class fails the Amalgamation Property \cite[Lemma~12]{GlassPierce}.
Their witness is a V-formation
\[
 G\hookrightarrow H,
 \qquad
 G\hookrightarrow K,
\]
in which all three groups have the same least strictly positive element $1$ and the embeddings preserve it.
In a discrete ordered group, the immediate successor and predecessor of $x$ are $x+1$ and $x-1$, respectively.
Hence these embeddings are positive normal embeddings.
Conversely, every positive normal embedding preserves $1$, since $1$ is the immediate successor of $0$.
Therefore an amalgam by positive normal embeddings would be an amalgam in the Glass--Pierce class, contrary to their lemma.
\end{proof}

\begin{theorem}\label{thm:kappaJ-failure}
Each of the classes
$$
 \Inv_{\odd}^{\Chain},
 \qquad
 \Inv_{\eveni}^{\Chain},
 \qquad
 \Inv_{\evenn}^{\Chain}
$$
fails the Amalgamation Property.
More generally, if
$$
 \Inv_{\varepsilon}^{\Chain}
 \subseteq
 \mathcal C
 \subseteq
 \Inv_{\oddeven}^{\Chain}
$$
for some $\varepsilon\in\{\odd,\eveni,\evenn\}$, then $\mathcal C$ fails AP.
\end{theorem}

\begin{proof}
Choose a non-amalgamable V-formation
$$
 \mathbf G_X\xrightarrow{\alpha}\mathbf G_Y,
 \qquad
 \mathbf G_X\xrightarrow{\beta}\mathbf G_Z
$$
as in Lemma~\ref{lem:discrete-group-obstruction}.

For the odd case, use the common two-element skeleton
$$
 \{t<u\},
 \qquad
 t\in\kappa_{\odd},
 \qquad
 u\in\kappa_J.
$$
Put the trivial group at $t$, put $\mathbf G_X,\mathbf G_Y,\mathbf G_Z$ at $u$, and take the transition from $t$ to $u$ to be trivial.
Use the identity maps on the skeletons and at the trivial layer, and use $\alpha,\beta$ at $u$.
The representation conditions are satisfied because the $u$-groups are discrete and the two layer maps are positive normal embeddings.
This gives a V-formation in $\Inv_{\odd}^{\Chain}$.

Suppose it had an amalgam in a class
$$
 \Inv_{\odd}^{\Chain}
 \subseteq
 \mathcal C
 \subseteq
 \Inv_{\oddeven}^{\Chain}.
$$
The target is necessarily odd, because embeddings preserve the constants $t$ and $f$.
Transport the amalgam to bunches.
The two copies of the common index $u$ have the same image $w$, and $w\in\kappa_J$ because bunch embeddings preserve the partition.
The $w$-layer is therefore a discrete abelian ordered group, and the two induced layer embeddings are positive and normal.
They amalgamate the group V-formation from Lemma~\ref{lem:discrete-group-obstruction}, a contradiction.

For the even idempotent-falsum case, use the same two-element skeleton with partition
$$
 \kappa_I=\{t\},
 \qquad
 \kappa_J=\{u\}.
$$
Put the trivial group and the trivial distinguished subgroup at $t$, place $\mathbf G_X,\mathbf G_Y,\mathbf G_Z$ at $u$, and take the transitions from $t$ to $u$ to be trivial.
Use $\alpha$ and $\beta$ on the $u$-layers.
The preceding layer argument at $u$ is unchanged.

For the even non-idempotent-falsum case, use the one-element skeleton
$$
 \{t\},
 \qquad
 t\in\kappa_J,
$$
and place the three discrete groups at its unique layer.
Any amalgam again yields a forbidden discrete ordered-group amalgam at the common $t$-layer.

In each even case, an amalgam in $\Inv_{\oddeven}^{\Chain}$ has the same even subtype, because the partition containing the least skeleton element is preserved by bunch embeddings.
This proves all assertions.
\end{proof}

\begin{remark}\label{rem:kappaJ}
The proof identifies the precise role of a $\kappa_J$-layer.
The obstruction is not merely that its group is linearly ordered.
It is the conjunction of discreteness and preservation of neighboring elements, which imports the category of discrete abelian linearly ordered groups with normal embeddings into the chain representation.
\end{remark}

\subsection{The distinguished-subgroup obstruction}

The absence of $\kappa_J$-layers does not imply AP.

\begin{theorem}\label{thm:symmetric-failure}
The classes
$$
 \Sym_{\odd}^{\Chain},
 \qquad
 \Sym_{\eveni}^{\Chain},
 \qquad\text{and}\qquad
 \Sym_{\odd}^{\Chain}\cup\Sym_{\eveni}^{\Chain}
$$
do not have the Amalgamation Property.
The counterexamples have exactly two positive idempotents.
\end{theorem}

\begin{proof}
We first construct an odd counterexample in the equivalent category of symmetric bunches.
Use the common skeleton
$$
 \{t<u\},
 \qquad
 t\in\kappa_{\odd},
 \qquad
 u\in\kappa_I.
$$
In each bunch let the $t$-layer be the trivial group, and let the transition from $t$ to $u$ be trivial.
At the $u$-layer take, with the usual order on $\Z$,
$$
 (G_u^{\mathcal X},H_u^{\mathcal X})=(\Z,\Z),
 \qquad
 (G_u^{\mathcal Y},H_u^{\mathcal Y})=(\Z,2\Z),
 \qquad
 (G_u^{\mathcal Z},H_u^{\mathcal Z})=(\Z,\Z).
$$
Define both maps out of $\mathcal X$ on the $u$-layer by
$$
 n\longmapsto2n,
$$
and use the identity maps on the skeleton and on the trivial $t$-layer.
These maps are symmetric-bunch embeddings.
They map the source subgroup into the respective target subgroup, and the complement condition is vacuous because
$$
 G_u^{\mathcal X}\setminus H_u^{\mathcal X}=\varnothing.
$$

Suppose that this V-formation has a bunch amalgam $\mathcal W$.
The two copies of $u$ have the same image $w\in\kappa_I^{\mathcal W}$.
Let
$$
 p:G_u^{\mathcal Y}\to G_w^{\mathcal W},
 \qquad
 q:G_u^{\mathcal Z}\to G_w^{\mathcal W}
$$
be the induced layer embeddings.
Commutativity gives
$$
 p(2)=q(2).
$$
Every linearly ordered abelian group is torsion-free, so
$$
 p(1)=q(1).
$$
On the other hand, the exact inverse-image condition gives
$$
 p(1)\notin H_w^{\mathcal W},
 \qquad
 q(1)\in H_w^{\mathcal W},
$$
because
$$
 1\notin2\Z=H_u^{\mathcal Y},
 \qquad
 1\in\Z=H_u^{\mathcal Z}.
$$
This is impossible.

For the even idempotent-falsum case, use the same construction with
$$
 \kappa_I=\{t,u\}
$$
and put
$$
 G_t=H_t=\{0\}
$$
in all three bunches.
The contradiction at the $u$-layer is unchanged.

Theorem~\ref{thm:representation} transfers these V-formations to the corresponding chain classes.
Each skeleton has two elements, hence each chain has exactly two positive idempotents.
Finally, an odd V-formation cannot acquire an even amalgam, or conversely, because bunch embeddings preserve the partition containing the least skeleton element.
Hence the union also fails AP.
\end{proof}

\subsection{Pure layers do not control a varying skeleton}

\begin{definition}\label{def:pure}
Let $A\leq B$ be torsion-free abelian groups.
The subgroup $A$ is \emph{pure} in $B$ if
$$
 A\cap mB=mA
 \qquad(m\geq1).
$$
Equivalently, $B/A$ is torsion-free.
An ordered-group embedding is \emph{pure} if its image is pure in the group reduct of its codomain.
A bunch embedding is \emph{layerwise pure} if every induced layer embedding is pure.
\end{definition}

\begin{proposition}\label{prop:pure-varying-failure}
There is a V-formation in $\Sym_{\odd}^{\Chain}$ whose two legs are layerwise pure and which has no amalgam in the full class
$$
 \Inv_{\oddeven}^{\Chain}
 =
 \Inv_{\odd}^{\Chain}
 \cup
 \Inv_{\eveni}^{\Chain}
 \cup
 \Inv_{\evenn}^{\Chain}.
$$
\end{proposition}

\begin{proof}
Use additive notation.
Let
$$
 A=\Z^2
$$
with the lexicographic order, and put
$$
 e_1=(1,0),
 \qquad
 e_2=(0,1).
$$

Define $\mathcal X$ on
$$
 \{t<s\},
$$
with
$$
 t\in\kappa_{\odd},
 \qquad
 s\in\kappa_I,
$$
by
$$
 G_t^{\mathcal X}=A,
 \qquad
 G_s^{\mathcal X}=H_s^{\mathcal X}=\{0\},
$$
and take the transition $t\to s$ to be zero.

Define $\mathcal Y$ on
$$
 \{t<p<s\},
$$
with
$$
 t\in\kappa_{\odd},
 \qquad
 p,s\in\kappa_I,
$$
by
$$
 G_t^{\mathcal Y}=G_p^{\mathcal Y}=A,
 \qquad
 H_p^{\mathcal Y}=2\Z\times\Z,
 \qquad
 G_s^{\mathcal Y}=H_s^{\mathcal Y}=\{0\}.
$$
Let the transition $t\to p$ be
$$
 a\longmapsto2a,
$$
and let every transition with target $s$ be zero.
This is a symmetric bunch because
$$
 2A\subseteq H_p^{\mathcal Y}.
$$

Define $\mathcal Z$ analogously on
$$
 \{t<q<s\},
$$
replacing $p$ by $q$ and putting
$$
 H_q^{\mathcal Z}=\Z\times2\Z.
$$

Embed $\mathcal X$ into $\mathcal Y$ and $\mathcal Z$ by the identity on the $t$-layer and the unique map on the trivial $s$-layer.
On skeletons use
$$
 t\mapsto t,
 \qquad
 s\mapsto s.
$$
Both maps are symmetric-bunch embeddings and are layerwise pure because their component maps are isomorphisms.

Suppose, toward a contradiction, that the corresponding chain V-formation has an amalgam
$$
 \mathbf W\in\Inv_{\oddeven}^{\Chain}.
$$
Since the two target chains are odd and embeddings preserve the constants $t$ and $f$, the amalgamating chain $\mathbf W$ must also be odd.
Transporting the amalgam through Theorem~\ref{thm:representation} yields an amalgamating bunch $\mathcal W$.
This bunch need not be symmetric; nevertheless, bunch embeddings preserve the skeleton partition, so the images of $p$ and $q$ still belong to $\kappa_I^{\mathcal W}$.
Let $p^*$ and $q^*$ be the images of $p$ and $q$ in its skeleton.
Write
$$
 \lambda_Y:A\to G_{p^*}^{\mathcal W},
 \qquad
 \lambda_Z:A\to G_{q^*}^{\mathcal W}
$$
for the layer embeddings, and write
$$
 \lambda_t:A\to G_t^{\mathcal W}
$$
for their common restriction at $t$.

If
$$
 p^*=q^*=r,
$$
compatibility with the two transitions gives
$$
 \lambda_Y(2a)
 =
 \sigma_{tr}^{\mathcal W}\lambda_t(a)
 =
 \lambda_Z(2a)
 \qquad(a\in A).
$$
Torsion-freeness gives
$$
 \lambda_Y(a)=\lambda_Z(a).
$$
But
$$
 e_1\notin H_p^{\mathcal Y},
 \qquad
 e_1\in H_q^{\mathcal Z},
$$
contradicting preservation of the subgroup and its complement.

If
$$
 p^*<q^*,
$$
the same calculation, followed by cancellation of $2$, gives
$$
 \sigma_{p^*q^*}^{\mathcal W}\lambda_Y(a)
 =
 \lambda_Z(a)
 \qquad(a\in A).
$$
Since $q^*\in\kappa_I^{\mathcal W}$, every transition into the $q^*$-layer has image in $H_{q^*}^{\mathcal W}$.
Hence
$$
 \lambda_Z(A)\subseteq H_{q^*}^{\mathcal W},
$$
contrary to complement preservation because
$$
 e_2\notin H_q^{\mathcal Z}.
$$
The case
$$
 q^*<p^*
$$
is symmetric, using
$$
 e_1\notin H_p^{\mathcal Y}.
$$
The skeleton is totally ordered, so these three cases are exhaustive and each yields a contradiction.
Hence the V-formation has no amalgam in $\Inv_{\oddeven}^{\Chain}$.
\end{proof}

\begin{remark}\label{rem:three-obstructions}
Theorems~\ref{thm:kappaJ-failure} and~\ref{thm:symmetric-failure}, and Proposition~\ref{prop:pure-varying-failure}, isolate different defects.
The first uses discrete normal layers.
The second uses distinguished subgroups on one fixed skeleton.
The third survives layerwise purity and is caused by incompatible enlargement of the skeleton.
\end{remark}

\subsection{Failure of strong amalgamation}

\begin{theorem}\label{thm:sap-failure}
The classes
$$
 \Sym_{\odd}^{\Chain},
 \qquad
 \Sym_{\eveni}^{\Chain},
 \qquad\text{and}\qquad
 \Sym_{\odd}^{\Chain}\cup\Sym_{\eveni}^{\Chain}
$$
fail the Strong Amalgamation Property.
More generally, every class
$$
 \mathcal C\subseteq\Inv_{\oddeven}^{\Chain}
$$
which contains either $\Sym_{\odd}^{\Chain}$ or $\Sym_{\eveni}^{\Chain}$ fails SAP.
\end{theorem}

\begin{proof}
The proof of \cite[Theorem~2]{CherriPowell} constructs a V-formation of totally ordered cyclic abelian groups having no strong amalgam even among representable lattice-ordered groups; hence it has no strong amalgam among abelian linearly ordered groups.
Choose a V-formation
$$
 \mathbf G_X\xrightarrow{\alpha}\mathbf G_Y,
 \qquad
 \mathbf G_X\xrightarrow{\beta}\mathbf G_Z
$$
which has no strong amalgam in that class.

For the odd case, form one-layer symmetric bunches on
$$
 \kappa=\{t\},
 \qquad
 t\in\kappa_{\odd},
$$
with layer groups
$$
 G_t^{\mathcal X}=\mathbf G_X,
 \qquad
 G_t^{\mathcal Y}=\mathbf G_Y,
 \qquad
 G_t^{\mathcal Z}=\mathbf G_Z,
$$
and layer embeddings $\alpha,\beta$.

For the even idempotent-falsum case, put
$$
 \kappa=\{t\},
 \qquad
 t\in\kappa_I,
$$
use the same three layer groups, and take
$$
 H_t^{\mathcal X}
 =
 H_t^{\mathcal Y}
 =
 H_t^{\mathcal Z}
 =
 \{0\}.
$$
The inverse-image condition follows from injectivity.

Suppose that either chain V-formation has a strong amalgam in a class
$$
 \mathcal C\subseteq\Inv_{\oddeven}^{\Chain}.
$$
The target has the same odd or even type because bunch embeddings preserve the partition containing the least skeleton element.
Its $t$-layer gives an amalgam of the original group V-formation.

This layer amalgam is strong.
Indeed, suppose
$$
 y\in G_t^{\mathcal Y},
 \qquad
 z\in G_t^{\mathcal Z},
 \qquad
 \mu_Y(y)=\mu_Z(z).
$$
The strong-amalgam condition gives an element $x$ of the common chain such that
$$
 y=\iota_Y(x),
 \qquad
 z=\iota_Z(x).
$$
In the odd one-layer case, the whole layer is the group $G_t$.
In the even $\kappa_I$-case, the embeddings preserve the distinction between dotted and undotted elements, while $y$ and $z$ are undotted.
Hence $x\in G_t^{\mathcal X}$ in either case.
Thus the $t$-layer amalgam is strong, contradicting the choice of the group V-formation.

The union fails SAP because a V-formation cannot mix the two types.
\end{proof}

\section{Pushouts and compatible orders}\label{sec:group-tools}

We use additive notation throughout this section.
Recall that a partially ordered abelian group is an abelian group equipped with a translation-invariant partial order.
Every linearly ordered abelian group is torsion-free.

\subsection{Pointwise pushouts, including non-injective input}

\begin{proposition}\label{prop:pointwise-pushout}
Let $\kappa$ be a directed poset and let
$$
 \mathcal X\xrightarrow{\mathbf f}\mathcal Y,
 \qquad
 \mathcal X\xrightarrow{\mathbf g}\mathcal Z
$$
be morphisms of direct systems of abelian groups over the same index set $\kappa$.
The component maps need not be injective.
For each $u\in\kappa$, form the abelian-group pushout
$$
 P_u=(Y_u\oplus Z_u)/D_u,
 \qquad
 D_u=\{(f_u(x),-g_u(x)):x\in X_u\},
$$
with canonical maps $\mu_{Y,u}$ and $\mu_{Z,u}$.
There are unique maps
$$
 \rho_{uv}:P_u\to P_v
$$
making
$$
\begin{aligned}
 \rho_{uv}\mu_{Y,u}
 &=
 \mu_{Y,v}\sigma^{\mathcal Y}_{uv},\\
 \rho_{uv}\mu_{Z,u}
 &=
 \mu_{Z,v}\sigma^{\mathcal Z}_{uv}.
\end{aligned}
$$
They make
$$
 (P_u,\rho_{uv})_{\kappa}
$$
a direct system and make the canonical families morphisms of direct systems.
Moreover,
$$
\begin{aligned}
 \mu_{Y,u}\text{ is injective}
 &\quad\Longleftrightarrow\quad
 \ker g_u\subseteq\ker f_u,\\
 \mu_{Z,u}\text{ is injective}
 &\quad\Longleftrightarrow\quad
 \ker f_u\subseteq\ker g_u.
\end{aligned}
$$
\end{proposition}

\begin{proof}
The two maps from $Y_u$ and $Z_u$ to $P_v$ obtained by first applying the transitions agree on $X_u$, by naturality of $\mathbf f$ and $\mathbf g$.
The pushout universal property therefore gives $\rho_{uv}$.
Uniqueness gives
$$
 \rho_{uu}=\mathrm{id},
 \qquad
 \rho_{uw}=\rho_{vw}\rho_{uv}.
$$

From the displayed quotient,
$$
 \ker\mu_{Y,u}=f_u(\ker g_u),
 \qquad
 \ker\mu_{Z,u}=g_u(\ker f_u).
$$
The two equivalences follow immediately.
\end{proof}

\begin{remark}\label{rem:transitions-not-enough}
Proposition~\ref{prop:pointwise-pushout} resolves the construction of the pointwise transition system for arbitrary input maps.
It does not by itself produce an amalgam of bunches.
The canonical legs may fail to be injective, the groups $P_u$ may have torsion, and the exact inverse-image condition for the distinguished subgroups may fail.
\end{remark}

\subsection{A simultaneous order-extension theorem}

\begin{theorem}\label{thm:simultaneous-order}
Let
$$
 (L_u,\sigma_{uv})_{u\leq v\text{ in }\kappa}
$$
be a direct system of torsion-free partially ordered abelian groups over an arbitrary chain $\kappa$, with positive transitions.
Each partial order can be extended to a total group order so that all the transitions remain positive.
The groups and the transition homomorphisms are unchanged.
\end{theorem}

\begin{proof}
For each $u\in\kappa$, let $P_u^0$ be the positive cone of the given partial order on $L_u$.
Let $\mathcal P$ be the set of all families
$$
 p=(P_u)_{u\in\kappa}
$$
such that $P_u$ is the positive cone of a partial group order on $L_u$ extending the given order, and
$$
 \sigma_{uv}(P_u)\subseteq P_v
 \qquad(u\leq v).
$$
The family $(P_u^0)_{u\in\kappa}$ belongs to $\mathcal P$, so $\mathcal P$ is nonempty.
Order $\mathcal P$ coordinatewise by inclusion.

Every chain $\mathcal C$ in $\mathcal P$ has the upper bound
$$
 \left(\bigcup_{p\in\mathcal C}P_u^p\right)_{u\in\kappa},
$$
where $P_u^p$ denotes the $u$-coordinate of $p$.
Indeed, the coordinatewise unions contain $0$, are closed under addition, and meet their negatives only in $\{0\}$: any two elements whose membership must be checked already occur together in one member of the chain.
The same observation shows that all transitions remain positive.
Zorn's lemma therefore gives a maximal family
$$
 (P_u)_{u\in\kappa}\in\mathcal P.
$$
Put $Q_u=P_u\setminus\{0\}$.
We prove that every $P_u$ is the positive cone of a total order.

Suppose that $P_s$ is not total.
Choose incomparable elements $a,b\in L_s$ and put $x=b-a$.
Then
$$
x\notin P_s
\qquad\text{and}\qquad
-x\notin P_s.
$$
Write $\mathbb N_0=\{0,1,2,\ldots\}$.
For $\varepsilon\in\{1,-1\}$ and $u\geq s$, set
$$
 a_u^\varepsilon=\sigma_{su}(\varepsilon x),
 \qquad
 P_u^\varepsilon
   =P_u+\mathbb N_0a_u^\varepsilon
   =\{p+na_u^\varepsilon:p\in P_u,\ n\in\mathbb N_0\},
$$
and, for $u<s$, set $P_u^\varepsilon=P_u$.
We claim that, for at least one common choice of $\varepsilon\in\{1,-1\}$, every $P_u^\varepsilon$ is a proper positive cone.

We first record the criterion used to prove the claim.
Let $P$ be the positive cone of a partially ordered torsion-free abelian group, put $Q=P\setminus\{0\}$, and let $a$ be an element of the group.
Then $P+\mathbb N_0a$ is a submonoid containing $P$, and
\[
 P+\mathbb N_0a
 \quad\text{fails to be a proper positive cone}
 \quad\Longleftrightarrow\quad
 -ma\in Q\text{ for some }m\geq1.
 \tag{*}\label{eq:cone-criterion}
\]
For the forward implication from right to left, if $-ma\in Q$, then both $ma$ and $-ma$ lie in $P+\mathbb N_0a$.
Conversely, suppose that a nonzero $c$ and its negative both lie in $P+\mathbb N_0a$.
Write
$$
 c=p+ra,
 \qquad
 -c=q+\ell a,
$$
where $p,q\in P$ and $r,\ell\in\mathbb N_0$.
Adding gives
$$
 -(r+\ell)a=p+q\in P.
$$
Here $r+\ell>0$, since otherwise $c,-c\in P$, contrary to the properness of $P$.
Moreover, $(r+\ell)a\neq0$: if it were zero, torsion-freeness would give $a=0$, and again $c,-c\in P$.
Hence $-(r+\ell)a\in Q$, proving~\eqref{eq:cone-criterion}.

Assume that neither sign works simultaneously at all later coordinates.
By~\eqref{eq:cone-criterion}, there are $u,v\geq s$ and $k,\ell\geq1$ such that
$$
 -k\sigma_{su}(x)\in Q_u,
 \qquad
 \ell\sigma_{sv}(x)\in Q_v.
$$
Since $\kappa$ is a chain, suppose first that $u\leq v$; the other case is symmetric.
Positivity of $\sigma_{uv}$ yields
$$
 -k\sigma_{sv}(x)\in P_v.
$$
This element is nonzero, since otherwise torsion-freeness would imply $\sigma_{sv}(x)=0$, contradicting $\ell\sigma_{sv}(x)\in Q_v$.
Thus $-k\sigma_{sv}(x)\in Q_v$.
Consequently both
$$
 k\ell\sigma_{sv}(x)
 \qquad\text{and}\qquad
 -k\ell\sigma_{sv}(x)
$$
are nonzero elements of $P_v$, a contradiction.
This proves the claim.

Fix a sign $\varepsilon$ supplied by the claim.
Each $P_u^\varepsilon$ contains $P_u$ and is a proper positive cone.
The transitions remain positive.
Indeed, if $s\leq r\leq u$, then
$$
 \sigma_{ru}(a_r^\varepsilon)=a_u^\varepsilon,
$$
and hence, for $p\in P_r$ and $n\in\mathbb N_0$,
$$
 \sigma_{ru}(p+na_r^\varepsilon)
 =\sigma_{ru}(p)+na_u^\varepsilon
 \in P_u^\varepsilon.
$$
If $r<s\leq u$, the cone at $r$ is unchanged and
$$
 \sigma_{ru}(P_r^\varepsilon)
 =\sigma_{ru}(P_r)
 \subseteq P_u
 \subseteq P_u^\varepsilon.
$$
If $r\leq u<s$, both cones are unchanged.
Thus $(P_u^\varepsilon)_{u\in\kappa}$ belongs to $\mathcal P$.

This extension is strict at $s$, because
$$
 \varepsilon x=a_s^\varepsilon\in P_s^\varepsilon
 \qquad\text{but}\qquad
 \varepsilon x\notin P_s.
$$
This contradicts maximality.
Therefore every $P_u$ defines a total group order.
By construction these orders extend the original ones, and all the original transitions remain positive.
\end{proof}

\begin{remark}[An alternative compactness proof]
The preceding proof is intrinsic to ordered algebra: it constructs a maximal compatible family of positive cones and shows directly that each of them is total.
For completeness, we give an alternative, shorter proof using the Compactness Theorem.
Although this second argument is model-theoretic rather than order-algebraic, it provides a useful independent formulation of the simultaneous order-extension principle.
\end{remark}

\begin{proof}[Alternative proof of Theorem~\ref{thm:simultaneous-order}]
We first treat one map.
Let
$$
 h:A\to B
$$
be a positive homomorphism of torsion-free partially ordered abelian groups, and suppose that the order on $B$ has been extended to a total group order.
The subgroup
$$
 K=\ker h
$$
is torsion-free, so its induced partial order extends to a total group order \cite[Corollary~13, p.~39]{Fuchs}.
Order $A$ lexicographically by
$$
 a<b
 \quad\Longleftrightarrow\quad
 h(a)<h(b)
 \quad\text{or}\quad
 \bigl(h(a)=h(b)\text{ and }0<b-a\text{ in }K\bigr).
$$
This is a total group order.
It extends the original order on $A$: if $a<b$ originally, either $h(a)<h(b)$ (by positivity of $h$) or $h(a)=h(b)$, in which case $b-a\in K$ is already $\geq 0$ for $K$'s order inherited from $A$, hence also $\geq 0$ for the extension of that order supplied by \cite[Corollary~13, p.~39]{Fuchs}.
Consequently $h$ is positive for the two extended orders.

Now let
$$
 F=\{u_1<\cdots<u_m\}
$$
be a finite subchain of $\kappa$.
Extend the order of $L_{u_m}$ to a total group order, again using \cite[Corollary~13, p.~39]{Fuchs}.
Working backwards, apply the one-map construction to
$$
 \sigma_{u_i u_{i+1}}:L_{u_i}\to L_{u_{i+1}}.
$$
Every transition between indices in $F$ is then positive by composition.

For arbitrary $\kappa$, use a many-sorted first-order language with one sort for each $u$, the group operations, the transitions, constants naming all elements of all $L_u$, and a relation $\leq_u$ on each sort.
Let $T$ contain all group-operation and transition equations true of the named elements, all inequations between constants naming distinct elements, the axioms saying that each $\leq_u$ is a total group order, all inequalities from the original partial orders, and the positivity axioms for every transition.
Every finite subset of $T$ involves only a finite subchain and is satisfiable by the preceding construction.
Compactness gives a model of $T$.
Restricting its orders to the named copies of the original groups gives the required total orders on the unchanged $L_u$'s.
\end{proof}

\subsection{Pure amalgamation over one skeleton}

\begin{lemma}\label{lem:same-skeleton}
Let
$$
 \mathcal X\xrightarrow{\mathbf f}\mathcal Y,
 \qquad
 \mathcal X\xrightarrow{\mathbf g}\mathcal Z
$$
be a V-formation of symmetric bunches.
Suppose that, after identifying the three partitioned skeletons by order isomorphisms, they have the same skeleton $\kappa$, both skeleton maps are the identity, and both bunch embeddings are layerwise pure.
Then the V-formation has an amalgamating symmetric bunch $\mathcal W$ over the same partitioned skeleton.
If every layer group of $\mathcal Y$ and $\mathcal Z$ is divisible, then every layer group of $\mathcal W$ is divisible.
\end{lemma}

\begin{proof}
Fix $u\in\kappa$, abbreviate the component embeddings by
$$
 f_u:X_u\to Y_u,
 \qquad
 g_u:X_u\to Z_u,
$$
and form their pushout in abelian groups:
$$
 P_u=(Y_u\oplus Z_u)/D_u,
 \qquad
 D_u=\{(f_u(x),-g_u(x)):x\in X_u\}.
$$
Write $\mu_{Y,u}$ and $\mu_{Z,u}$ for the canonical maps.

The group $P_u$ is torsion-free.
Indeed, $D_u$ is pure in $Y_u\oplus Z_u$.
Let $m\geq1$ and suppose
$$
 m(y,z)=(f_u(x),-g_u(x)).
$$
Purity of $f_u(X_u)$ in $Y_u$ gives $x_0\in X_u$ such that
$$
 f_u(x)=m f_u(x_0).
$$
Since $f_u$ is injective, $x=mx_0$.
Torsion-freeness of $Y_u$ and $Z_u$ then gives
$$
 (y,z)=(f_u(x_0),-g_u(x_0))\in D_u.
$$
Thus
$$
 (Y_u\oplus Z_u)/D_u
$$
is torsion-free.
Since $f_u$ and $g_u$ are injective, so are the canonical maps.
The quotient normal form also gives
\begin{equation}\label{eq:pushout-intersection}
 \mu_{Y,u}(Y_u)\cap\mu_{Z,u}(Z_u)
 =
 \mu_{Y,u}f_u(X_u)
 =
 \mu_{Z,u}g_u(X_u).
\end{equation}

Let $Y_u^+$ and $Z_u^+$ be the positive cones of the two factor orders, and put
$$
 C_u
 =
 \mu_{Y,u}(Y_u^+)
 +
 \mu_{Z,u}(Z_u^+).
$$
This is a proper positive cone.
Suppose that
$$
 w\in C_u\cap(-C_u).
$$
Then
$$
 w=\mu_{Y,u}(y_1)+\mu_{Z,u}(z_1)
$$
and
$$
 -w=\mu_{Y,u}(y_2)+\mu_{Z,u}(z_2)
$$
for nonnegative $y_1,y_2,z_1,z_2$.
Hence
$$
 \mu_{Y,u}(y_1+y_2)
 =
 -\mu_{Z,u}(z_1+z_2).
$$
By \eqref{eq:pushout-intersection}, there is $x\in X_u$ such that
$$
 f_u(x)=y_1+y_2,
 \qquad
 g_u(x)=-(z_1+z_2).
$$
The first equality gives $x\geq0$, and the second gives $x\leq0$, because both component maps are order embeddings.
Thus $x=0$.
Hence $y_1+y_2=z_1+z_2=0$; since all four elements are nonnegative, they are all zero, and consequently $w=0$.

The cone $C_u$ induces exactly the original order on each factor.
For example, suppose
$$
 \mu_{Y,u}(y)\in C_u.
$$
Then
$$
 \mu_{Y,u}(y)
 =
 \mu_{Y,u}(y_+)+\mu_{Z,u}(z_+)
$$
for some
$$
 y_+\geq0,
 \qquad
 z_+\geq0.
$$
By \eqref{eq:pushout-intersection}, there is $x\in X_u$ such that
$$
 y-y_+=f_u(x),
 \qquad
 z_+=g_u(x).
$$
The second equality gives $x\geq0$, and therefore
$$
 y\geq y_+\geq0.
$$
The argument for the $Z$-factor is symmetric.
Thus $P_u$, ordered by $C_u$, is a torsion-free partially ordered abelian group into which both factors order-embed.

For $u\leq v$, Proposition~\ref{prop:pointwise-pushout} supplies a unique transition
$$
 \rho_{uv}:P_u\to P_v.
$$
It is positive for the cones $C_u$, because it is positive on both generating factors.
We have therefore obtained a direct system of torsion-free partially ordered abelian groups.
Apply Theorem~\ref{thm:simultaneous-order} to extend all orders simultaneously to total group orders.
The canonical maps remain ordered-group embeddings, because their images were already totally ordered by the partial orders.

For $u\in\kappa_I$, define
\begin{equation}\label{eq:pushout-H}
 H_u^{\mathcal W}
 =
 \mu_{Y,u}(H_u^{\mathcal Y})
 +
 \mu_{Z,u}(H_u^{\mathcal Z}).
\end{equation}
If $s<u$, then $\rho_{su}$ maps each factor of $P_s$ into the corresponding subgroup on the right of \eqref{eq:pushout-H}.
Hence
$$
 \rho_{su}(P_s)\subseteq H_u^{\mathcal W}.
$$

It remains to check complements.
Suppose
$$
 \mu_{Y,u}(y)\in H_u^{\mathcal W}.
$$
Then
$$
 \mu_{Y,u}(y-h_Y)
 =
 \mu_{Z,u}(h_Z)
$$
for some
$$
 h_Y\in H_u^{\mathcal Y},
 \qquad
 h_Z\in H_u^{\mathcal Z}.
$$
By \eqref{eq:pushout-intersection}, there is $x\in X_u$ such that
$$
 y-h_Y=f_u(x),
 \qquad
 h_Z=g_u(x).
$$
The inverse-image condition for the original $Z$-leg gives
$$
 x\in H_u^{\mathcal X}.
$$
The same condition for the $Y$-leg then gives
$$
 y\in H_u^{\mathcal Y}.
$$
Therefore
$$
 \mu_{Y,u}^{-1}(H_u^{\mathcal W})
 =
 H_u^{\mathcal Y}.
$$
The proof for the $Z$-leg is symmetric.

After applying the preceding tagging convention, the groups $P_u$, transitions $\rho_{uv}$, partitioned skeleton, and subgroups \eqref{eq:pushout-H} form a symmetric bunch $\mathcal W$.
The canonical families are bunch embeddings and agree on $\mathcal X$.

Finally, if $Y_u$ and $Z_u$ are divisible, then $Y_u\oplus Z_u$ is divisible, and so is its quotient $P_u$.
The order-extension step does not change the group reducts.
\end{proof}

\section{The fixed-idempotent positive theorem}\label{sec:positive}

\begin{definition}\label{def:positive-class}
An idempotent-symmetric odd or even involutive $\FLe$-chain $\mathbf X$ is \emph{layer-divisible} if every canonical group
$$
 G_u(\mathbf X),
 \qquad
 u\in\Idp(\mathbf X),
$$
is divisible as an abelian group.
For $n\geq1$ and
$$
 \varepsilon\in\{\odd,\eveni\},
$$
put
$$
 \Div_{\varepsilon,n}^{\Chain}
 =
 \left\{
 \mathbf X\in\Sym_{\varepsilon}^{\Chain}:
 |\Idp(\mathbf X)|=n
 \text{ and $\mathbf X$ is layer-divisible}
 \right\},
$$
and set
$$
 \Div_{\odd\eveni,n}^{\Chain}
 =
 \Div_{\odd,n}^{\Chain}
 \cup
 \Div_{\eveni,n}^{\Chain}.
$$
\end{definition}

\begin{remark}\label{rem:intrinsic}
Definition~\ref{def:positive-class} imposes two properties on each algebra.
It does not require a chosen V-formation to have the same skeleton or pure legs.
Those representation-side facts follow automatically in the proof.
\end{remark}

\begin{theorem}\label{thm:main}
For every finite $n\geq1$, each of the classes
$$
 \Div_{\odd,n}^{\Chain},
 \qquad
 \Div_{\eveni,n}^{\Chain},
 \qquad\text{and}\qquad
 \Div_{\odd\eveni,n}^{\Chain}
$$
has the Amalgamation Property.
\end{theorem}

\begin{proof}
Fix
$$
 \varepsilon\in\{\odd,\eveni\}
$$
and a V-formation in
$$
 \Div_{\varepsilon,n}^{\Chain}.
$$
By Theorem~\ref{thm:representation}, it becomes a V-formation
$$
 \mathcal X\xrightarrow{\mathbf f}\mathcal Y,
 \qquad
 \mathcal X\xrightarrow{\mathbf g}\mathcal Z
$$
of symmetric bunches.
Each skeleton has exactly $n$ elements.

The skeleton maps of bunch embeddings are injective order maps.
An injective map between two finite sets of the same cardinality is bijective.
Hence both skeleton maps are order isomorphisms.
After transporting the data along these isomorphisms, the three bunches have the same partitioned skeleton and both skeleton maps are the identity.

Fix a skeleton point $u$.
The source group $X_u$ is divisible.
Its image under either layer embedding is therefore a divisible subgroup of the torsion-free target group.
Such a subgroup is pure.
Indeed, let $m\geq1$ and suppose that $f_u(x)\in mY_u$, say
$$
 f_u(x)=my.
$$
Choose $x_0\in X_u$ with
$$
 x=mx_0.
$$
Then
$$
 m\bigl(y-f_u(x_0)\bigr)=0,
$$
and torsion-freeness gives
$$
 y=f_u(x_0).
$$
The same argument applies to $g_u$.
Thus both bunch embeddings are layerwise pure.

Lemma~\ref{lem:same-skeleton} supplies an amalgamating symmetric bunch $\mathcal W$ over the same $n$-element skeleton.
Every layer group of $\mathcal W$ is divisible because the layer groups of $\mathcal Y$ and $\mathcal Z$ are divisible.
Theorem~\ref{thm:representation} transports $\mathcal W$ and its two embeddings back to an amalgam in
$$
 \Div_{\varepsilon,n}^{\Chain}.
$$

For the union, an embedding preserves the partition containing the least skeleton element.
Hence there is no embedding between an odd member and an even member.
Every V-formation in
$$
 \Div_{\odd\eveni,n}^{\Chain}
$$
lies entirely in one of its two summands, where the preceding argument applies.
\end{proof}

\begin{remark}\label{rem:roles}
The hypotheses in Theorem~\ref{thm:main} have distinct roles.
Fixed finite $n$ converts the induced skeleton embeddings into isomorphisms, excluding the mechanism of Proposition~\ref{prop:pure-varying-failure}.
Layer divisibility makes the source images pure and the pushout layers divisible.
Neither hypothesis is merely a technical condition on the arrows of one chosen diagram.
\end{remark}

\section{Semilinear variety-level consequences}\label{sec:varieties}

An extension
$$
 \mathbf A\leq\mathbf B
$$
is \emph{essential} if
$$
 \theta\cap A^2\neq\Delta_A
 \qquad
 \text{for every $\theta\in\operatorname{Con}(\mathbf B)$ with
 $\theta\neq\Delta_B$.}
$$
Equivalently, every homomorphism from $\mathbf B$ whose restriction to $\mathbf A$ is injective is itself injective.
An embedding is essential if its image is an essential subalgebra.
A V-formation
$$
 \mathbf A\longrightarrow\mathbf B,
 \qquad
 \mathbf A\longrightarrow\mathbf C
$$
is \emph{essential} if the designated second leg $\mathbf A\to\mathbf C$ is essential.
A class has the \emph{Essential Amalgamation Property} (EAP) if every essential V-formation in the class has an amalgam in the class.

\begin{lemma}\label{lem:semilinear-transfer}
Let $\mathcal V$ be a variety of semilinear commutative pointed residuated lattices.
Then $\mathcal V$ has AP if and only if the class of its finitely subdirectly irreducible members, equivalently its linearly ordered members, has EAP.
\end{lemma}

\begin{proof}
This is \cite[Corollary~3.5]{FussnerSantschi}; it is based on the general transfer results of \cite{FussnerMetcalfe}.
\end{proof}

\begin{lemma}[Finite-rank bunch data for finitely generated chains]
\label{lem:finitely-generated-chain-bunch}
Let $\mathbf X$ be a finitely generated odd or even involutive $\FLe$-chain.
Then its canonical bunch has a finite skeleton, every layer group $G_u$ is a finitely generated free abelian group, and every distinguished subgroup $H_u$, $u\in\kappa_I$, is finitely generated \cite[\S18]{FuchsInfinite}.
\end{lemma}

\begin{proof}
Let $S$ be a finite generating set of $\mathbf X$, and put
$$
\tau(x)=x\to x.
$$
Since $f=t'$ and $\tau(t)=t$, \cite[Lemma~3.1(ix)]{JeneiRep} gives
$$
\tau(f)=t.
$$
Let $x\in X$, and choose a term over $S$ whose value is $x$.
After fixing this evaluation, recursively replace each occurrence of $\wedge$ or $\vee$ by the argument selected by that lattice operation.
The resulting term has the same value and contains only multiplication, residuals, generators, and constants.
By \cite[Lemma~3.1(xvi)]{JeneiRep}, $\tau(x)$ is the maximum of the $\tau$-values of the generator and constant occurrences in this term.
Consequently,
$$
\kappa(\mathbf X)
=
\{\tau(x):x\in X\}
\subseteq
\{\tau(s):s\in S\}\cup\{t\},
$$
so the skeleton is finite.

For $v\leq u$ in the skeleton, let
$$
\rho^{\mathbf X}_{v\to u}:X_v\longrightarrow X_u
$$
denote the canonical layer transition.
For each $u\in\kappa(\mathbf X)$, put
$$
E_u
=
\bigl\{
 \rho^{\mathbf X}_{\tau(s)\to u}(s):
 s\in S,\ \tau(s)\leq u
\bigr\}
\cup
\bigl\{
 \rho^{\mathbf X}_{t\to u}(t),
 \rho^{\mathbf X}_{t\to u}(f)
\bigr\}.
$$
We claim that the finite set $E_u$ generates the $u$-layer algebra.

Indeed, let $x\in X_u$ and prune a term representing $x$ as above.
Since $\tau(x)=u$, \cite[Lemma~3.1(xvi)]{JeneiRep} shows that every generator or constant occurring in the pruned term belongs to a layer indexed by some $v\leq u$.
The coherence of the canonical transitions and the layer-reconstruction formulas in \cite[Lemma~4.2]{JeneiRep} show, by induction on the term, that its value is unchanged if every occurrence from a $v$-layer is first mapped into $X_u$ by $\rho^{\mathbf X}_{v\to u}$ and all operations are then evaluated in $X_u$.
Thus $x$ belongs to the layer subalgebra generated by $E_u$, proving the claim.

Suppose first that $u\in\kappa_{\odd}\cup\kappa_J$.
By \cite[Lemmas~3.2, 5.2, and~7.2]{JeneiRep}, the canonical group $G_u$ has the same universe and multiplication as the corresponding cancellative layer algebra; at a $\kappa_J$-point only the designated falsum is changed when passing to the associated odd algebra.
In its ordered-group reduct,
$$
x^{-1}=x\to u,
\qquad
x\to y=x^{-1}y.
$$
Since the order is total, $\wedge$ and $\vee$ select one of their arguments.
Induction on layer terms therefore shows that $E_u$ generates $G_u$ as an abelian group.

Now suppose that $u\in\kappa_I$.
The canonical homomorphism of the $u$-layer maps it surjectively onto the cancellative odd algebra $\pi_1(X_u)$, whose ordered-group reduct is $G_u$ \cite[Theorem~6.3 and Lemma~7.2]{JeneiRep}.
Since $E_u$ generates $X_u$, its image generates $\pi_1(X_u)$.
The same term-to-group argument therefore shows that $G_u$ is finitely generated in this case as well.

Every $G_u$ is a finitely generated torsion-free abelian group and hence is free abelian of finite rank \cite[\S15]{FuchsInfinite}.
For $u\in\kappa_I$, the subgroup $H_u\leq G_u$ is finitely generated because every subgroup of a finitely generated abelian group is finitely generated.
\end{proof}

\begin{remark}[Finite chains versus finite-rank chain data]
\label{rem:finite-versus-finite-rank}
For the odd/even failures considered here, the chain reduction in Lemma~\ref{lem:semilinear-transfer} cannot be restricted to finite chains.
Indeed, every canonical layer group of a finite odd or even involutive $\FLe$-chain is finite and hence trivial, since every abelian $o$-group is torsion-free.
Moreover, the trivial ordered group is not discrete, so the corresponding skeleton has no $\kappa_J$-points.
It follows from \cite[Example~8.2]{JeneiRep} that every such finite chain is a finite Sugihara chain.

The varieties of Sugihara monoids and of odd Sugihara monoids have AP \cite[Theorem~4.5 and Corollary~4.6]{FussnerSantschiSemilinear}.
Hence, by Lemma~\ref{lem:semilinear-transfer}, their chain members have EAP.
Consequently, an essential V-formation all of whose members are finite odd chains amalgamates among odd Sugihara chains, while one all of whose members are finite odd or even chains amalgamates among Sugihara chains.
Such V-formations therefore cannot witness failure of EAP in the full odd or combined odd/even chain regimes.
This excludes only V-formations lying entirely in those regimes; finite essential witnesses involving other involutive chains occur in \cite[Theorem~5.2]{FussnerSantschiSemilinear}.

Corollary~3.7 of \cite{FussnerSantschiSemilinear} reduces AP instead to essential V-formations of finitely generated chains.
By Lemma~\ref{lem:finitely-generated-chain-bunch}, these chains have finite-skeleton bunches with finitely generated free-abelian layer groups.
Neither result supplies a uniform bound on the skeleton sizes or layer ranks, so this reduction does not by itself yield a finite effective test for AP.
\end{remark}
\bigskip

Write
$$
 \Inv_{\odd}^{\Sem}
$$
for the variety of odd semilinear involutive commutative $\FLe$-algebras.
For a class $\mathcal K$ of algebras, write $\mathsf V(\mathcal K)$ for the variety generated by $\mathcal K$.

\begin{lemma}\label{lem:essential-counterexamples}
The following non-amalgamable V-formations can be chosen essential on the leg into $\mathbf Z$:
\begin{enumerate}[(i)]
\item the odd discrete-layer V-formation used in Theorem~\ref{thm:kappaJ-failure};
\item the odd and even distinguished-subgroup V-formations used in Theorem~\ref{thm:symmetric-failure};
\item an even one-layer distinguished-subgroup V-formation.
\end{enumerate}
None of these V-formations has an amalgam in any odd or even involutive $\FLe$-chain.
\end{lemma}

\begin{proof}

We shall repeatedly use the embedding characterization \cite[Lemma~10]{JeneiDense}.
In particular, if the skeleton map of a bunch homomorphism is an order embedding, all its layer components are ordered-group embeddings, and the skeleton partition, the distinguished subgroups and their complements, and the neighbourhood operations at $\kappa_J$-layers are preserved, then the corresponding $\FLe$-algebra homomorphism is an embedding.

To check essentiality, it is enough to consider surjective homomorphisms onto chain images.
Indeed, given a homomorphism $\psi:\mathbf Z\to\mathbf W$, replace its codomain by the subalgebra $\psi[Z]$.
This image is again a chain.
A homomorphic image of an odd chain is odd.
A homomorphic image of an even chain is odd if the images of $f$ and $t$ coincide, and is even otherwise; in the latter case idempotence of $f$ is preserved.
To verify the lower-cover condition, suppose that
\[
 \psi(f)<\psi(x)<\psi(t).
\]
Since the source is a chain, monotonicity forces $f<x<t$, contradicting that $f$ is the lower cover of $t$.
For (i), use the odd two-layer V-formation from the proof of Theorem~\ref{thm:kappaJ-failure} and choose
$$
 \iota_Z:\mathbf X\to\mathbf Z.
$$
Let
$$
 \psi:\mathbf Z\to\mathbf W
$$
be a homomorphism into an involutive $\FLe$-chain, and suppose that $\psi$ is not injective.
Because $\mathbf Z$ is odd and homomorphisms preserve both constants, $\mathbf W$ is odd.
Pass to bunches.
The skeleton of $\mathbf Z$ is
$$
 \{t<u\},
 \qquad
 t\in\kappa_{\odd},
 \qquad
 u\in\kappa_J.
$$
The image of $u$ under the skeleton map is either in $\kappa_J^{\mathcal W}$ or in $\kappa_{\odd}^{\mathcal W}$.

If the image of $u$ belongs to $\kappa_J^{\mathcal W}$, then the induced homomorphism on the discrete $u$-layer preserves neighboring elements.
It is injective on that layer.
Indeed, its kernel is a convex subgroup \cite[Chapter~II, \S5]{Fuchs}.
If the kernel contained a nonzero element, then it would contain the least strictly positive element of the discrete source group.
But preservation of upper neighbors maps that element to the least strictly positive element of the target group, rather than to the identity, a contradiction.
It is also injective on the trivial $t$-layer.
If $\alpha(t)=\alpha(u)$, then, since $t$ and $u$ are idempotent elements equal to their own skeleton points already lying in the image of $\mathbf X$, this gives $\psi(t)=\psi(u)$ directly, so the restriction of $\psi$ to $\mathbf X$ is not injective, proving essentiality of $\iota_Z$.
Otherwise the two skeleton points have distinct images, so the underlying direct-system map is an embedding.
The skeleton partition and the neighbourhood operations are preserved, so \cite[Lemma~10]{JeneiDense} shows that $\psi$ is an embedding, contrary to assumption.

Therefore the image of $u$ belongs to $\kappa_{\odd}^{\mathcal W}$.
The image of $t$ also belongs to $\kappa_{\odd}^{\mathcal W}$, and that part contains at most the least skeleton point.
Thus
$$
 \psi(u)=\psi(t).
$$
The two distinct skeleton units $u$ and $t$ already belong to the common subalgebra $\mathbf X$.
Hence the restriction of $\psi$ to $\mathbf X$ is not injective.
This proves essentiality of $\iota_Z$.

For (ii), consider first the odd distinguished-subgroup V-formation from Theorem~\ref{thm:symmetric-failure}, again with the leg into $\mathbf Z$ distinguished.
Its skeleton is
$$
 \{t<u\},
 \qquad
 t\in\kappa_{\odd},
 \qquad
 u\in\kappa_I,
$$
and
$$
 H_u^{\mathcal Z}=G_u^{\mathcal Z}=\Z.
$$
If a bunch homomorphism out of $\mathcal Z$ maps $u$ to a $\kappa_I$-point, the partial-injectivity condition for bunch homomorphisms, applied inside $H_u^{\mathcal Z}=G_u^{\mathcal Z}$, makes its $u$-component strictly order-preserving and hence injective.
The trivial $t$-component is also injective.
If $\alpha(t)=\alpha(u)$, then, since $t$ and $u$ are idempotent elements equal to their own skeleton points already lying in the image of $\mathbf X$, this gives $\psi(t)=\psi(u)$ directly, so the restriction of the corresponding algebra homomorphism to $\mathbf X$ is not injective.
Otherwise the two skeleton points have distinct images, so the underlying direct-system map is an embedding; preservation of the skeleton partition and of the distinguished subgroups and their complements allows us to apply \cite[Lemma~10]{JeneiDense}, and the corresponding algebra homomorphism is an embedding.
If instead $u$ is mapped to the unique $\kappa_{\odd}$-point, then $u$ and $t$ are identified; these two skeleton units already lie in the image of $\mathcal X$.
Thus every non-injective homomorphism out of $\mathbf Z$ is non-injective on the image of $\mathbf X$, and the leg is essential.
For the even two-layer construction, corestrict the homomorphism to its chain image.
If this image is odd, then the homomorphism identifies the dotted and undotted copies of $0$ in the common $t$-layer, so its restriction to the common algebra is already non-injective.

Suppose therefore that the image is even.
Then $t$ is mapped to its least $\kappa_I$-point.
If $u$ is mapped to a distinct $\kappa_I$-point, the $u$-component is strictly order-preserving and the trivial $t$-layer map is injective.
The skeleton map is then an order embedding preserving the partition, and the distinguished subgroups and their complements are preserved.
Hence \cite[Lemma~10]{JeneiDense} shows that the whole homomorphism is an embedding.
If $t$ and $u$ have the same image, then the restriction to the common algebra identifies its two skeleton units.
These are the only possibilities, since an even bunch has no odd skeleton point.

For (iii), use the one-point even skeleton
$$
 \{t\},
 \qquad t\in\kappa_I,
$$
and put
$$
 (G_t^{\mathcal X},H_t^{\mathcal X})=(\Z,\Z),
 \qquad
 (G_t^{\mathcal Y},H_t^{\mathcal Y})=(\Z,2\Z),
 \qquad
 (G_t^{\mathcal Z},H_t^{\mathcal Z})=(\Z,\Z),
$$
with both embeddings out of $\mathcal X$ given by $m\mapsto2m$.
The calculation in the proof of Theorem~\ref{thm:symmetric-failure} applies verbatim: in any alleged amalgam, cancellation gives equality of the two images of $1$, while the exact inverse-image condition places that element outside and inside the amalgamating distinguished subgroup.
The leg into $\mathcal Z$ is essential as well.
If a homomorphism maps $t$ to a $\kappa_I$-point, its layer component is strictly order-preserving because $H_t^{\mathcal Z}=G_t^{\mathcal Z}$.
It is therefore an ordered-group embedding, and preservation of the distinguished subgroup and its complement allows \cite[Lemma~10]{JeneiDense} to be applied.
Hence the corresponding algebra homomorphism is an embedding.
If it maps $t$ to the unique $\kappa_{\odd}$-point, it identifies the two constants of the even common algebra and is therefore non-injective already there.

Finally, an odd chain cannot embed into an even chain, or conversely, because embeddings preserve the constants and their equality or inequality.
The non-amalgamation calculations above therefore rule out amalgams in the full class of odd or even involutive chains, not merely in the symmetric subclasses.
\end{proof}

\begin{theorem}\label{thm:essential-span-transfer}
Let $\mathcal V$ be a variety of semilinear involutive commutative pointed residuated lattices.
\begin{enumerate}[(i)]
\item
If $\mathcal V$ contains the three chains of one of the odd V-formations in Lemma~\ref{lem:essential-counterexamples}, then $\mathcal V$ fails AP.

\item
The same conclusion holds for one of the even V-formations provided every linearly ordered member of $\mathcal V$ is odd or even.
\end{enumerate}
\end{theorem}

\begin{proof}
A chain receiving an embedding from an odd algebra is itself odd.
Under the hypothesis in \textup{(ii)}, a chain receiving an embedding from an even algebra cannot be odd and hence is even.
Lemma~\ref{lem:essential-counterexamples} therefore rules out a chain amalgam in either case, and Lemma~\ref{lem:semilinear-transfer} applies.
\end{proof}

\begin{corollary}\label{cor:symmetric-variety-failure}
Each of the varieties
$$
 \mathsf V(\Sym_{\odd}^{\Chain}),
 \qquad
 \mathsf V(\Sym_{\eveni}^{\Chain}),
 \qquad\text{and}\qquad
 \mathsf V(\Sym_{\odd}^{\Chain}\cup\Sym_{\eveni}^{\Chain})
$$
fails the Amalgamation Property.
More generally, any semilinear variety containing the odd distinguished-subgroup V-formation fails AP.
A semilinear variety containing an even distinguished-subgroup V-formation fails AP whenever all its linearly ordered members are odd or even.
\end{corollary}

\begin{proof}
Each displayed variety is semilinear, since semilinearity is equational for commutative pointed residuated lattices \cite{FussnerSantschiSemilinear} and every generator is a chain.
The odd-generated variety and the variety generated by the union contain the odd essential V-formation, so Theorem~\ref{thm:essential-span-transfer}\textup{(i)} applies.

For $\mathsf V(\Sym_{\eveni}^{\Chain})$, every linearly ordered member is odd or even.
Indeed, commutative pointed residuated lattices are congruence-distributive, and the finitely subdirectly irreducible members of a semilinear variety are exactly its chains \cite[Section~2]{FussnerSantschiSemilinear}.
Hence the FSI form of Jónsson's Lemma \cite[Lemma~1.1]{AglianoMarcos} shows that every such chain is a homomorphic image of a subalgebra of an ultraproduct of the generating even chains.
Ultraproducts and subalgebras remain even, while a homomorphic image of an even chain is either even or odd according as the two constants remain distinct or become identified.
Theorem~\ref{thm:essential-span-transfer}\textup{(ii)} now applies.
\end{proof}

\begin{corollary}\label{cor:semilinear-failure}
Let $\mathcal V$ be a variety of semilinear involutive commutative pointed residuated lattices.
If $\mathcal V$ contains the odd discrete-layer essential V-formation from Lemma~\ref{lem:essential-counterexamples}, then $\mathcal V$ fails the Amalgamation Property.
Consequently, every such variety satisfying
$$
\Inv_{\odd}^{\Sem}\subseteq\mathcal V
$$
fails the Amalgamation Property.

In particular, the variety of semilinear involutive commutative residuated lattices fails the Amalgamation Property, recovering the recent result of Fussner and Santschi \cite[Theorem~5.2]{FussnerSantschiSemilinear}.
\end{corollary}

\begin{proof}
The first assertion follows directly from Theorem~\ref{thm:essential-span-transfer}\textup{(i)}.
The odd discrete-layer V-formation belongs to $\Inv_{\odd}^{\Sem}$, so every variety containing $\Inv_{\odd}^{\Sem}$ contains this obstructing V-formation.
The final assertion is obtained by taking $\mathcal V$ to be the variety of all semilinear involutive commutative residuated lattices.
\end{proof}

The preceding results also determine what happens after closing the positive classes of Section~\ref{sec:positive} under the operations of a variety.

\begin{theorem}\label{thm:generated-divisible-varieties}
The varieties generated by the fixed-$n$, layer-divisible chain classes have the following AP status:
\begin{enumerate}[(i)]
\item $\mathsf V(\Div_{\odd,1}^{\Chain})$ has AP;
\item $\mathsf V(\Div_{\odd,n}^{\Chain})$ fails AP for every $n\geq2$;
\item $\mathsf V(\Div_{\eveni,n}^{\Chain})$ and $\mathsf V(\Div_{\odd\eveni,n}^{\Chain})$ fail AP for every $n\geq1$.
\end{enumerate}
\end{theorem}

\begin{proof}
Each displayed variety is semilinear, since semilinearity is equational for commutative pointed residuated lattices and every generator is a chain.
An odd one-layer chain is term-equivalent to a linearly ordered abelian group, and layer divisibility is precisely divisibility of that group.
Every linearly ordered abelian group embeds, as an ordered subgroup, into a divisible linearly ordered abelian group, while every abelian lattice-ordered group is a subdirect product of linearly ordered abelian groups.
Consequently, $\mathsf V(\Div_{\odd,1}^{\Chain})$ is term-equivalent to the variety of abelian lattice-ordered groups, which has AP by \cite[Corollaries~6.2 and~6.3]{FussnerSantschiSemilinear}.

For $n\geq2$, each of the three chains in the odd distinguished-subgroup V-formation embeds into a member of $\Div_{\odd,n}^{\Chain}$.
Indeed, in each corresponding bunch, replace the nontrivial layer group $\Z$ by $\mathbb Q$ and retain the respective distinguished subgroup $\Z$, $2\Z$, or $\Z$ as a subgroup of $\mathbb Q$.
The inclusion $\Z\hookrightarrow\mathbb Q$ satisfies the required exact inverse-image condition for these distinguished subgroups.
When $n>2$, adjoin $n-2$ trivial divisible $\kappa_I$-layers above the original nontrivial $\kappa_I$-layer, with trivial distinguished subgroups and zero transitions.
The resulting enlarged bunch belongs to $\Div_{\odd,n}^{\Chain}$, and the original bunch embeds into it.
Hence each of the three original chains is a subalgebra of a member of $\Div_{\odd,n}^{\Chain}$.
Since varieties are closed under subalgebras, all three chains belong to $\mathsf V(\Div_{\odd,n}^{\Chain})$.
Theorem~\ref{thm:essential-span-transfer}\textup{(i)} therefore gives failure of AP.

For the even case, use the one-layer V-formation from Lemma~\ref{lem:essential-counterexamples}.
In each corresponding bunch, replace the layer group $\Z$ by $\mathbb Q$ and retain the respective distinguished subgroup $\Z$, $2\Z$, or $\Z$ as a subgroup of $\mathbb Q$.
The inclusion $\Z\hookrightarrow\mathbb Q$ satisfies the required exact inverse-image condition for these distinguished subgroups.
When $n>1$, adjoin $n-1$ trivial divisible $\kappa_I$-layers above the original layer, with trivial distinguished subgroups and zero transitions.
The resulting enlarged bunch belongs to $\Div_{\eveni,n}^{\Chain}$, and the original one-layer bunch embeds into it.
Hence all three chains are subalgebras of members of $\Div_{\eveni,n}^{\Chain}$.
Since varieties are closed under subalgebras, the same essential V-formation belongs both to $\mathsf V(\Div_{\eveni,n}^{\Chain})$ and to $\mathsf V(\Div_{\odd\eveni,n}^{\Chain})$.
Every linearly ordered member of either generated variety is odd or even.
For the even-generated variety this follows exactly as in Corollary~\ref{cor:symmetric-variety-failure}.
For the variety generated by the union, let $\mathbf D$ be a linearly ordered member.
By the same FSI form of Jónsson's Lemma \cite[Lemma~1.1]{AglianoMarcos}, $\mathbf D$ is a homomorphic image of a subalgebra of an ultraproduct
$$
 \prod_{i\in I}\mathbf A_i/\mathcal U,
$$
where every $\mathbf A_i$ belongs to the odd or the even generating class.
The sets of odd and even indices partition $I$, so exactly one of them belongs to $\mathcal U$.
Restricting to that $\mathcal U$-large set shows that the ultraproduct is odd or even.
These respective properties pass to subalgebras.
A homomorphic image of an odd chain is odd, while a homomorphic image of an even chain is even or odd according as the two constants remain distinct or become identified.
Hence Theorem~\ref{thm:essential-span-transfer}\textup{(ii)} shows that both varieties fail AP.
\end{proof}

\begin{remark}\label{rem:class-variety-contrast}
Theorem~\ref{thm:main} and Theorem~\ref{thm:generated-divisible-varieties} are not in conflict.
The former concerns V-formations and amalgams inside a non-equational class of chains having exactly $n$ positive idempotents and divisible layers.
Variety generation introduces subalgebras, homomorphic images, and products; except in the odd one-layer case, this recovers an essential non-amalgamable span.

All the varieties shown above to fail AP also fail SAP, since SAP implies AP.
They also fail the Transferable Injections Property: for varieties of commutative pointed residuated lattices, AP and that property are equivalent by \cite[Corollary~44]{MetcalfeMontagnaTsinakis}.
\end{remark}

\section{Conclusion}\label{sec:conclusion}

The layer-group representation reveals a hierarchy of amalgamation obstructions.

A $\kappa_J$-layer imports the discrete ordered-group obstruction through normal embeddings.
Eliminating $\kappa_J$ removes that obstruction but leaves the distinguished-subgroup obstruction at $\kappa_I$-layers.
Even layerwise purity does not repair arbitrary diagrams, because the targets may enlarge the skeleton incompatibly.

The positive theorem addresses the latter two difficulties by intrinsic restrictions on the algebras.
A fixed finite number of positive idempotents makes all skeleton embeddings in a V-formation isomorphisms.
Divisible canonical layers make the induced layer embeddings pure and ensure that the pushout layers remain divisible.
The simultaneous order-extension theorem then assembles the layer pushouts into an amalgamating bunch.

The negative and positive results therefore separate three structural obstruction mechanisms from a positive regime.
The unrestricted full classes fail AP because of discrete layers.
The unrestricted idempotent-symmetric classes fail AP because of distinguished subgroups, whereas the fixed-$n$, layer-divisible idempotent-symmetric subclasses have AP.
Strong amalgamation still fails in the unrestricted idempotent-symmetric chain classes.

At the variety level, Theorem~\ref{thm:essential-span-transfer} packages the chain counterexamples into a general transfer principle: a semilinear variety containing one of the odd essential non-amalgamable spans constructed here fails AP, and the analogous conclusion holds for an even span whenever all of its linearly ordered members are odd or even.
This yields, in particular, failure of AP for the varieties generated by the idempotent-symmetric odd chains, the idempotent-symmetric even chains, and their union, as well as for every semilinear variety containing $\Inv_{\odd}^{\Sem}$.
Moreover, variety generation preserves the positive result only in the odd one-layer case.
For every $n\geq2$ the variety generated by $\Div_{\odd,n}^{\Chain}$ fails AP, and for every $n\geq1$ the varieties generated by $\Div_{\eveni,n}^{\Chain}$ and $\Div_{\odd\eveni,n}^{\Chain}$ fail AP.
The exceptional variety $\mathsf V(\Div_{\odd,1}^{\Chain})$ is term-equivalent to the variety of abelian lattice-ordered groups and has AP.
All of the negative variety cases consequently fail SAP and, by the AP--TIP equivalence for varieties of commutative pointed residuated lattices, the Transferable Injections Property.

\section*{Acknowledgements}
This work was supported by the Ministry of Culture and Innovation of Hungary, through the National Research, Development and Innovation Fund, grant no.~K138596.

\end{document}